\documentclass[11pt]{amsart}

\usepackage{color, amsmath,amssymb, amsfonts, amstext,amsthm, latexsym}

   \usepackage{latexsym, oldlfont}

   \newtheorem{lemma}{Lemma}[section]
   \newtheorem{theorem}[lemma]{Theorem}

   \newtheorem{coro}[lemma]{Corollary}
   \newtheorem{definition}[lemma]{Definition}

\newcommand{\R}{{\mathbb R}}

\newcommand{\beq}{\begin{equation}}
\newcommand{\eeq}{\end{equation}}
\newcommand{\Del}{\Delta}

\newcommand{\nab}{\nabla}
\newcommand{\lam}{\lambda}
\newcommand{\sig}{\sigma}

\newcommand{\ome}{\omega}
\newcommand{\patl}{\partial}

\newcommand{\al}{\alpha}

\title[Euler and Navier-Stokes Equations]{Non-uniqueness of Regular
Solutions for Incompressible Static Euler Equations with Given
Boundary Conditions and Turbulent Global Solutions of
Incompressible Navier-Stokes Equations}

\author{Yongqian Han}

\date{}

\address{
Institute of Applied Physics and Computational Mathematics,
Beijing 100088, China}
\address{
National Key Laboratory of Computational Physics,
Beijing 100088, China}
\email{han\_yongqian@iapcm.ac.cn}

\keywords{Incompressible Navier-Stokes Equations, Incompressible
Static Euler Equations, Randomness and Turbulence of Incompressible
Fluids, Typing Models without Prandtl Layer.}

\subjclass[2000]{35Q30, 76D05, 76F02, 37L20}

\begin{document}

\begin{abstract} The incompressible Navier-Stokes equations and static
Euler equations are considered. We find that there exist infinite
non-trivial regular solutions of incompressible static Euler equations 
with given boundary conditions. Moreover there exist random solutions 
of incompressible static Euler equations. Provided Reynolds number is
large enough and time variable $t$ goes to infinity, these random
solutions of static Euler equations are the path limits of
corresponding Navier-Stokes flows. But the double limits of these
Navier-Stokes flows do not exist. These phenomena reveal randomness
and turbulence of incompressible fluids. Therefore these solutions are
called turbulent solutions. Here some typing models without Prandtl
layer are given.
\end{abstract}

\maketitle

~ ~ ~ ~

\section{Introduction}
\setcounter{equation}{0}

The incompressible Navier-Stokes equations in $\R^3$ are given by
\beq\label{NS1}
u_t-\nu\Del u+(u\cdot\nab)u+\nab P=0,
\eeq
\beq\label{NS2}
\nab\cdot u=0,
\eeq
where $u=u(t,x)=\big(u^1(t,x),u^2(t,x),u^3(t,x)\big)$ and $P=P(t,x)$
stand for the unknown velocity vector field of fluid and its
pressure, $\nu>0$ is the coefficient of viscosity. Here $x=(x_1,x_2,
x_3)$, $\nab=(\patl_1,\patl_2,\patl_3)$ and $\patl_j=\frac{\patl}
{\patl x_j}\;\;$  ($j=1,2,3$).

Taking curl with equation \eqref{NS1}, we have
\beq\label{NS-G-ome}
\ome_t-\nu\Del\ome+(u\cdot\nab)\ome-(\ome\cdot\nab)u=0.
\eeq
Here vorticity $\ome=(\ome^1,\ome^2,\ome^3)$ and
\beq\label{Def-G-ome}\begin{split}
\ome(t,x)=&\nab\times u(t,x).\\
\end{split}\eeq

The incompressible static Euler equations in $\R^3$ are given by
\beq\label{SEE1}
(u\cdot\nab)u+\nab P=0,
\eeq
\beq\label{SEE2}
\nab\cdot u=0.
\eeq
By taking curl, the equation \eqref{SEE1} is equivalent to
\beq\label{SEE-G-ome}
(u\cdot\nab)\ome-(\ome\cdot\nab)u=0.
\eeq

The velocity vector $u$ is solution of equations \eqref{NS1}
\eqref{NS2} provided $(u,\ome)$ satisfies equation \eqref{NS-G-ome}.
The velocity vector $u$ is solution of equations \eqref{SEE1}
\eqref{SEE2} provided $(u,\ome)$ satisfies equation \eqref{SEE-G-ome}.

The Fourier transformation of $u(t,x)$ with respect to $x$ denotes by
$\hat{u}(t,\xi)$. Then equation \eqref{NS2} and \eqref{SEE2} are
equivalent to
\beq\label{NS2-F}
\xi_j\hat{u}^j=\xi_1\hat{u}^1+\xi_2\hat{u}^2+\xi_3\hat{u}^3=0.
\eeq
It is that $\hat{u}(t,\xi)$ is perpendicular to $\xi\in \R^3$. Denote
by $\xi\bot\hat{u}$ for any $\xi\in\R^3$.

There is a large literature (for instance, see \cite{Con01, YGuo23,
han1, han3, LiW19, Te01, QZha24} and therein references) studying the
incompressible Navier-Stokes and steady Euler equations.

To construct solutions of Navier-Stokes equations \eqref{NS1} \eqref{NS2}
and Euler equations \eqref{SEE1} \eqref{SEE2}, we define some linear
operators which are perpendicular to operator $\nab$.

\begin{definition}\label{AOp-def} Let domain $D\subset\R^3$. Assume
that a vector $v(x)=\big(v_1(x),v_2(x),v_3(x)\big)$ satisfies
equation \eqref{NS2} for any $x\in D$ and the Fourier transformation
$\hat{v}(\xi)$ is well defined.

First provided $D$ is bounded, we define a linear operator $v(\nab)$
as follows
\beq\begin{split}\label{AOp-def-f}
\big(v(\nab)f\big)\hat{}\,(\xi)=\hat{v}(\xi)\hat{f}(\xi),
\;\;\forall f\in C(D),
\end{split}\eeq
where $\hat{f}(\xi)$ is the Fourier transformation of $f(x)$.

Next provided $D=\R^3$, the linear operator $v(\nab)$ can be defined
in space $C^{\infty}_c(\R^3)$
\beq\begin{split}\label{AOp-def-ff}
\big(v(\nab)f\big)\hat{}\,(\xi)=\hat{v}(\xi)\hat{f}(\xi),\;\;
\forall f\in C^{\infty}_c(\R^3).
\end{split}\eeq

Especially let a vector $P=P(x)=\big(p_1(x),p_2(x),p_3(x)\big)$ and
$p_j(x)$ $(j=1,2,3)$ be polynomial. Provided polynomial vector $P$
satisfies
\beq\label{NS2-FP}
x\cdot P(x)=x_jp_j(x)=x_1p_1(x)+x_2p_2(x)+x_3p_3(x)\equiv 0,
\;\;\;\;\;\;\forall x\in\R^3.
\eeq
Then we define linear operator $P(\nab)=\big(p_1(\nab),p_2(\nab),
p_3(\nab)\big)$. Now we have
\beq\begin{split}\label{AOp-def-ffP}
&\big(P(\nab)f\big)\hat{}\,(\xi)=P(i\xi)\hat{f}(\xi)\\
=&\big(p_1(i\xi)\hat{f}(\xi),p_2(i\xi)\hat{f}(\xi),
p_3(i\xi)\hat{f}(\xi)\big),\;\;\forall f\in C^{\infty}_c(\R^3),
\end{split}\eeq
where $\hat{f}(\xi)$ is the Fourier transformation of $f(x)$.
\end{definition}

It is obvious that any linear operator $A(\nab)$ defined in Definition
\ref{AOp-def} is perpendicular to operator $\nab$. Let velocity vector
$u(x)$ be the solution of static Euler equations \eqref{SEE1} \eqref{SEE2}
and regular enough, then we can define a linear operator $u(\nab)$ by
using Definition \ref{AOp-def}.

Choose a real linear operator $A=A(\nab)=\big(a_1(\nab),a_2(\nab),
a_3(\nab)\big)$ be defined in Definition \ref{AOp-def}, then $A(\nab)$
is perpendicular to $\nab$, i.e. ~ $\nab\cdot A(\nab)=0$.
Since $\nab,A(\nab)$ and $A(\nab)\times\nab$ are vertical two by two,
we can write velocity vector $u$ as follows
\beq\begin{split}\label{Sym-Rep-a2}
u(t,x)&=A(\nab)\phi(t,x)+\{A(\nab)\times\nab\}\psi(t,x)\\
&=\Big(a_1(\nab)\phi+\big(a_2(\nab)\patl_3-a_3(\nab)\patl_2\big)\psi,
a_2(\nab)\phi+\big(a_3(\nab)\patl_1-a_1(\nab)\patl_3\big)\psi,\\
&\hspace{8mm}a_3(\nab)\phi+\big(a_1(\nab)\patl_2-a_2(\nab)\patl_1\big)\psi
\Big)(t,x),
\end{split}\eeq
where $\phi$ and $\psi$ are real scalar functions. Then $u$ satisfies
the equation \eqref{NS2} and \eqref{SEE2}. We call that the
formulation \eqref{Sym-Rep-a2} is symplectic representation
of velocity vector $u$.

Now the curl of velocity vector $u$ can be rewritten as follows
\beq\label{Def-SRa2-ome}\begin{split}
\ome(t,x)=&\nab\times u(t,x)\\
=&-\{(A(\nab)\times\nab\}\phi(t,x)-\{\big(A(\nab)\times\nab\big)\times\nab\}\psi(t,x)\\
=&-\{(A(\nab)\times\nab\}\phi(t,x)+A(\nab)\Del\psi(t,x)\\
=&\Big(-\big(a_2(\nab)\patl_3-a_3(\nab)\patl_2\big)\phi+a_1(\nab)\Del\psi,\\
&\hspace{4mm}-\big(a_3(\nab)\patl_1-a_1(\nab)\patl_3\big)\phi+a_2(\nab)\Del\psi,\\
&\hspace{4mm}-\big(a_1(\nab)\patl_2-a_2(\nab)\patl_1\big)\phi
+a_3(\nab)\Del\psi\Big)(t,x).\\
\end{split}\eeq

Some turbulent global solutions of equations \eqref{NS1} \eqref{NS2}
have been constructed in \cite{han2} by using \eqref{Sym-Rep-a2}
\eqref{Def-SRa2-ome} with special linear operator
$A(\nab)=e_3\times\nab$, provided $(\phi,\psi)$ is radial symmetric
and cylindrical symmetric respectively.

By employing equations \eqref{NS-G-ome} \eqref{SEE-G-ome}
\eqref{Sym-Rep-a2} \eqref{Def-SRa2-ome}, some new interesting results
are obtained in this paper. The main result is as follows

\begin{theorem}[Global Turbulent Solution]
\label{SRa2-Thm-GTS}  Assume that $A(\nab)$ is any linear operator
defined by Definition \ref{AOp-def}, $\lam$ is constant, $\psi(x)$ is
regular enough and satisfies the following equation
\beq\label{Eig-ps}\begin{split}
-\Del\psi(x)=\lam^2\psi(x),
\end{split}\eeq
and $\phi(x)=\lam\psi(x)$.

Let
\beq\begin{split}\label{Def-ue}
u_e(x)&=A(\nab)\phi(x)+\{A(\nab)\times\nab\}\psi(x)\\
&=\lam A(\nab)\psi(x)+\{A(\nab)\times\nab\}\psi(x),
\end{split}\eeq
then $u=u_e(x)$ is the solution of static incompressible Euler
equations \eqref{SEE1} \eqref{SEE2}.

Moreover let
\beq\begin{split}\label{Def-uns}
u_{ns}(t,x)&=e^{-\nu\lam^2t}A(\nab)\phi(x)
+e^{-\nu\lam^2t}\{A(\nab)\times\nab\}\psi(x)\\
&=\lam e^{-\nu\lam^2t}A(\nab)\psi(x)
+e^{-\nu\lam^2t}\{A(\nab)\times\nab\}\psi(x),\\
\end{split}\eeq
then $u=u_{ns}$ is the solution of incompressible Navier-Stokes
equations \eqref{NS1} \eqref{NS2}.
\end{theorem}

\ \ \

It is easy to prove this Theorem \ref{SRa2-Thm-GTS}. In fact,
inserting $\phi(x)=\lam\psi(x)$ and equation \eqref{Eig-ps} into
\eqref{Sym-Rep-a2} \eqref{Def-SRa2-ome}, we find $\ome=-\lam u$.
Inserting $\ome=-\lam u$ into equations \eqref{NS-G-ome}
\eqref{SEE-G-ome}, we can prove Theorem \ref{SRa2-Thm-GTS}.

\ \ \

By taking samples of $\big(A(\nab),\lam,\phi,\psi\big)$ in
Theorem \ref{SRa2-Thm-GTS}, some new interesting results, such as
non-uniqueness of static incompressible Euler equations \eqref{SEE1}
\eqref{SEE2} with given boundary conditions, can be shown clearly.

Let us solve the equation \eqref{Eig-ps} with periodic boundary
condition
\beq\label{ps-pb}\begin{split}
\psi(x+2\pi)=\psi(x),\;\;\;\;\;\;\forall x\in\R^3.
\end{split}\eeq
This is eigenvalue problem. By spectrum theory of compact operator
$(-\Del)^{-1}$, there exist countable infinite eigenvalues
\[
0\le\lam_1\le\cdots\le\lam_n\le\cdots,\;\;n=1,2,\cdots.
\]
Each eigenvalue is real, and may be geometrically simple or finite.

For eigenvalue $\lam_n$ of problem \eqref{Eig-ps} \eqref{ps-pb},
we define linear spaces $E_n$ and $E_n^{\ast}$ as follows
\beq\label{ps-pb-EF}\begin{split}
E_n=\big\{\psi_n(x)\big|&\mbox{any eigenfunction $\psi_n$
corresponding to eigenvalue $\lam_n$} \big\},
\end{split}\eeq
\beq\label{ps-pb-Ea}\begin{split}
E_n^{\ast}=\mbox{span}\{&\sin(k_1x_1)\sin(k_2x_2)\sin(k_3x_3),
\sin(k_1x_1)\sin(k_2x_2)\cos(k_3x_3),\\
&\sin(k_1x_1)\cos(k_2x_2)\sin(k_3x_3),
\sin(k_1x_1)\cos(k_2x_2)\cos(k_3x_3),\\
&\cos(k_1x_1)\sin(k_2x_2)\sin(k_3x_3),
\cos(k_1x_1)\sin(k_2x_2)\cos(k_3x_3),\\
&\cos(k_1x_1)\cos(k_2x_2)\sin(k_3x_3),
\cos(k_1x_1)\cos(k_2x_2)\cos(k_3x_3),\\
&\lam_n^2=k_1^2+k_2^2+k_3^2,\mbox{ integer } k_j\ge0,\;\;j=1,2,3\}.
\end{split}\eeq
It is obvious that $E_n^{\ast}\subseteq E_n$.

As the consequence of Theorem \ref{SRa2-Thm-GTS}, we have

\begin{theorem}[Global Turbulent Periodic Solution]
\label{SRa2P-Thm-GTPS}  Let $A(\nab)$ be any linear operator which is
defined in Definition \ref{AOp-def} by using polynomials. Provided
that $\lam_n$ is eigenvalue and $\psi_n$ is eigenfunction of problem
\eqref{Eig-ps} \eqref{ps-pb}, $\psi_n\in E_n$, $\phi_n(x)=\lam_n
\psi_n(x)$ $(n=1,2,\cdots)$. Let
\beq\begin{split}\label{Def-ue-Pn}
u_e^n(x)&=A(\nab)\phi_n(x)+\{A(\nab)\times\nab\}\psi_n(x)\\
&=\lam_n A(\nab)\psi_n(x)+\{A(\nab)\times\nab\}\psi_n(x),
\end{split}\eeq
then $u=u_e^n(x)$ $(n=1,2,\cdots)$ is the solution of static
incompressible Euler equations \eqref{SEE1} \eqref{SEE2} with periodic
boundary condition
\beq\label{VVuE-pb}\begin{split}
u(x+2\pi)=u(x),\;\;\forall x\in\R^3.
\end{split}\eeq

Moreover let
\beq\begin{split}\label{Def-uns-Pn}
u_{ns}^n(t,x)&=e^{-\nu\lam^2_nt}A(\nab)\phi_n(x)
+e^{-\nu\lam^2_nt}\{A(\nab)\times\nab\}\psi_n(x)\\
&=\lam_ne^{-\nu\lam^2_nt}A(\nab)\psi_n(x)
+e^{-\nu\lam^2_nt}\{A(\nab)\times\nab\}\psi_n(x),\\
\end{split}\eeq
then $u=u_{ns}^n(t,x)$ $(n=1,2,\cdots)$ is the solution of
incompressible Navier-Stokes equations \eqref{NS1} \eqref{NS2} with
periodic boundary condition
\beq\label{VVuNS-pb}\begin{split}
u(t,x+2\pi)=u(t,x),\;\;\forall t\ge0,\;x\in\R^3
\end{split}\eeq
and initial data
\beq\label{VVu-id}\begin{split}
u(t,x){\big|}_{t=0}=u_e^n(x).
\end{split}\eeq
\end{theorem}

It is obvious that the following corollaries are the consequences of
Theorem \ref{SRa2P-Thm-GTPS}.

\begin{coro}[Non-unique Periodic Solutions of Euler Equations]
\label{Nonu-SEE-pb} The solution of static incompressible Euler
equations \eqref{SEE1} \eqref{SEE2} with periodic boundary condition
\eqref{VVuE-pb} is non-unique.
\end{coro}

\begin{coro}[Global Well-posedness of NS Periodic Solution]
\label{GWP-NS-pb} The incompressible Navier-Stokes equations
\eqref{NS1} \eqref{NS2} with periodic boundary condition
\eqref{VVuNS-pb} and initial data
\beq\label{VVu-u0}\begin{split}
u(t,x){\big|}_{t=0}=u_0(x)
\end{split}\eeq
is global well-posed, provided $u_0\in B$. Here subset $B$ of
$C^{\infty}_{per}$ is defined by
\beq\label{SubS-B-L2}\begin{split}
B=\bigcup_{A(\nab)}\bigcup_{n=1,2,\cdots}\big\{u(x)\big|&
u=\lam_n A(\nab)\psi_n(x)+\{A(\nab)\times\nab\}\psi_n(x),
\forall\psi_n\in E_n\big\}.
\end{split}\eeq
$A(\nab)$ is any linear operator which is defined in Definition
\ref{AOp-def} by using polynomials.
\end{coro}

To show randomness and turbulence of incompressible fluids, we
consider limits of solution $u_{ns}^n$ of Navier-Stokes equations as
$\nu\rightarrow0$ and $t\rightarrow\infty$.

Let $\mu$ be random number and $\Gamma_{\mu}$ be a path. In the path
$\Gamma_{\mu}$, we have
\beq\label{Ran-Pa}\begin{split}
t\nu=\mu\\
\end{split}\eeq
for any $(\nu,t)\in\Gamma_{\mu}$.

For the solution $u_{ns}^n(t,x)$ of Navier-Stokes equations
constructed in \eqref{Def-uns-Pn}, we have path limit
\beq\label{PL-NSs-pb}\begin{split}
u_e^{n,\mu}(x)=&\lim_{(\nu,t)\in\Gamma_{\mu},(\nu,1/t)\rightarrow(0,0)}
u_{ns}^n(t,x)\\
=&\lam_ne^{-\mu\lam^2_n}A(\nab)\psi_n(x)
+e^{-\mu\lam^2_n}\{A(\nab)\times\nab\}\psi_n(x)\\
=&\lam_nA(\nab)\{e^{-\mu\lam^2_n}\psi_n(x)\}
+\{A(\nab)\times\nab\}\{e^{-\mu\lam^2_n}\psi_n(x)\}.\\
\end{split}\eeq
It is obvious that random $u_e^{n,\mu}$ is the solution of static
incompressible Euler equations \eqref{SEE1} \eqref{SEE2} with periodic
boundary condition \eqref{VVuE-pb} since $E_n$ is linear space and
$e^{-\mu\lam^2_n}\psi_n(x)\in E_n$. On the other hand double limit
$\lim_{(\nu,1/t)\rightarrow(0,0)}u_{ns}^n(t,x)$ does not exist. These
phenomena reveal randomness and turbulence of incompressible fluids.

\ \ \

Now we give another typification to show that Prandtl layer dose not
appear.

For any bounded domain $D\subset\R^3$ and $\patl D\in C^{\infty}$, we
consider equation \eqref{Eig-ps} with Dirichlet boundary condition
\beq\label{ps-Db}\begin{split}
\psi(x)\big|_{x\in\patl D}=0.
\end{split}\eeq
The equation \eqref{Eig-ps} with Dirichlet boundary condition
\eqref{ps-Db} is so-called  Dirichlet eigenvalue problem. The
eigenvalue $\lam$ of problem \eqref{Eig-ps} \eqref{ps-Db} is
so-called Dirichlet eigenvalue. By spectrum theory of compact
operator $(-\Del)^{-1}$, there exist countable infinite Dirichlet
eigenvalues
\[
0<\lam_1\le\cdots\le\lam_n\le\cdots,\;\; n=1,2,\cdots
\]
Each Dirichlet eigenvalue is real, and may be geometrically simple or
finite. We define linear space $F_n$ as follows
\beq\label{Ef-Db-n}\begin{split}
F_n=\big\{\psi_n(x)\big|\mbox{any eigenfunction $\psi_n$ corresponding
to eigenvalue $\lam_n$}\big\}.
\end{split}\eeq
By the regular theory of elliptic equations, we have $\psi_n\in
H^{\infty}(D)$ for any $\psi_n\in F_n$ and $n=1,2,\cdots$

As corollary of Theorem \ref{SRa2-Thm-GTS}, we derive

\begin{theorem}[Global Turbulent Dirichlet Solution]
\label{SRa2P-Thm-GTDS}  Let $A(\nab)$ be any linear operator which is
defined in Definition \ref{AOp-def}. Provided that $\lam_n$ is
eigenvalue and $\psi_n$ is eigenfunction of problem
\eqref{Eig-ps} \eqref{ps-Db}, $\psi_n\in F_n$, $\phi_n(x)=\lam_n
\psi_n(x)$ $(n=1,2,\cdots)$. Let
\beq\begin{split}\label{Def-ue-Dn}
u_e^{D,n}(x)&=A(\nab)\phi_n(x)+\{A(\nab)\times\nab\}\psi_n(x)\\
&=\lam_n A(\nab)\psi_n(x)+\{A(\nab)\times\nab\}\psi_n(x),
\end{split}\eeq
then $u=u_e^{D,n}(x)$ $(n=1,2,\cdots)$ is the solution of static
incompressible Euler equations \eqref{SEE1} \eqref{SEE2} with
boundary condition
\beq\label{EVVu-Db}\begin{split}
\sig\cdot u(x)\big|_{x\in\patl D}=\{\lam_n\sig\cdot A(\nab)\psi_n(x)
+\sig\cdot\{A(\nab)\times\nab\}\psi_n(x)\}\big|_{x\in\patl D},
\end{split}\eeq
where $\sig$ is the outer normal of $\patl D$.

Moreover let
\beq\begin{split}\label{Def-uns-Dn}
u_{ns}^{D,n}(t,x)&=e^{-\nu\lam^2_nt}A(\nab)\phi_n(x)
+e^{-\nu\lam^2_nt}\{A(\nab)\times\nab\}\psi_n(x)\\
&=\lam_ne^{-\nu\lam^2_nt}A(\nab)\psi_n(x)
+e^{-\nu\lam^2_nt}\{A(\nab)\times\nab\}\psi_n(x),\\
\end{split}\eeq
then $u=u_{ns}^{D,n}(t,x)$ $(n=1,2,\cdots)$ is the solution of
incompressible Navier-Stokes equations \eqref{NS1} \eqref{NS2} with
Dirichlet boundary condition
\beq\label{NSVVu-Db}\begin{split}
u(t,x)\big|_{x\in\patl D}=\{\lam_ne^{-\nu\lam^2_nt}A(\nab)\psi_n(x)
+e^{-\nu\lam^2_nt}\{A(\nab)\times\nab\}\psi_n(x)\}\big|_{x\in\patl D}
\end{split}\eeq
and initial data
\beq\label{NSVVu-Did}\begin{split}
u(t,x){\big|}_{t=0}=u_e^{D,n}(x).
\end{split}\eeq

Indeed the solution of Navier-Stokes equations \eqref{NS1} \eqref{NS2}
with Dirichlet boundary condition \eqref{NSVVu-Db} and initial data
\eqref{NSVVu-Did} converges to the solution of Euler equations
\eqref{SEE1} \eqref{SEE2} with boundary condition \eqref{EVVu-Db} as
$\nu\rightarrow0$, i.e.
\[
\lim_{\nu\rightarrow0}u_{ns}^{D,n}=u_e^{D,n}.
\]
Here Prandtl layer dose not appear.
\end{theorem}

Similarly we have path limit
\beq\label{PL-NSs-pb}\begin{split}
u_e^{D,n,\mu}(x)=&\lim_{(\nu,t)\in\Gamma_{\mu},(\nu,1/t)\rightarrow(0,0)}
u_{ns}^{D,n}(t,x)\\
=&\lam_ne^{-\mu\lam^2_n}A(\nab)\psi_n(x)
+e^{-\mu\lam^2_n}\{A(\nab)\times\nab\}\psi_n(x)\\
=&\lam_nA(\nab)\{e^{-\mu\lam^2_n}\psi_n(x)\}
+\{A(\nab)\times\nab\}\{e^{-\mu\lam^2_n}\psi_n(x)\}.\\
\end{split}\eeq
Since $F_n$ is linear space and $e^{-\mu\lam^2_n}\psi_n(x)\in F_n$,
it is obvious that random $u_e^{D,n,\mu}$ is the solution of static
incompressible Euler equations \eqref{SEE1} \eqref{SEE2} with
boundary condition
\beq\label{EVVu-Db1}\begin{split}
&\sig\cdot u(x)\big|_{x\in\patl D}\\
=&\big(\lam_n\sig\cdot A(\nab)\{e^{-\mu\lam^2_n}\psi_n(x)\}
+\sig\cdot\{A(\nab)\times\nab\}\{e^{-\mu\lam^2_n}\psi_n(x)\}
\big)\big|_{x\in\patl D}.
\end{split}\eeq
On the other hand double limit $\lim_{(\nu,1/t)\rightarrow(0,0)}
u_{ns}^{D,n}(t,x)$ does not exist. These phenomena also reveal
randomness and turbulence of incompressible fluids.

\ \ \

Assume that $A(\nab)=e_3\times\nab$, $\phi=0$ and $\psi(t,x)=
\psi(t,r,x_3)$ ($r^2=x_1^2+x_2^2$) in \eqref{Sym-Rep-a2}. Inserting
this special \eqref{Sym-Rep-a2} into Navier-Stokes equations
\eqref{NS1} \eqref{NS2}, the equations reduce to axially symmetric
Navier-Stokes equations without swirl. The global well-posedness of
axially symmetric Navier-Stokes equations without swirl was first
established by Ladyzhenskaya \cite{Lad68} and Ukhovskii-Yudovich
\cite{UYu68}. Leonardi et al. \cite{Leo99} give a simpler proof to
this result by using a priori weight estimates. Recently a different
proof of the global well-posed result was given in \cite{han3}.

\ \ \

Provided that $A(\nab)=e_3\times\nab=(-\patl_2,\patl_1,0)$,
$\phi(t,x)=\phi(x)$ and $\psi(t,x)=0$ in \eqref{Sym-Rep-a2}. Inserting
this special \eqref{Sym-Rep-a2} and \eqref{Def-SRa2-ome} into static
Euler equation \eqref{SEE-G-ome}, we derive

\begin{theorem}[Non-unique Solutions of Periodic Pipeline Euler Flow]
\label{Nonu-SEE-ppDb}  For any function $\phi(\xi,\eta)\in C^{\infty}
\big(\R\times(L_1,L_2)\big)$ with
\beq\label{ph-ppl}\begin{split}
\phi(\xi+2\pi,\eta)=\phi(\xi,\eta),\;\;\forall\xi\in\R,
\;\eta\in(L_1,L_2),
\end{split}\eeq
we define
\beq\begin{split}\label{Def-ue-ppl}
u_e^{pp}(x)&=\big(-\patl_2\phi(x_1+x_2,x_3),
\patl_1\phi(x_1+x_2,x_3),0\big),
\end{split}\eeq
then $u=u_e^{pp}(x)$ is the solution of static incompressible Euler
equations \eqref{SEE1} \eqref{SEE2} with boundary condition
\beq\label{ue-ppl-b}\begin{split}
&u(x_1+2\pi,x_2,x_3)=u(x_1,x_2,x_3),\;\forall(x_1,x_2)\in\R^2,
\;x_3\in(L_1,L_2),\\
&u(x_1,x_2+2\pi,x_3)=u(x_1,x_2,x_3),\;\forall(x_1,x_2)\in\R^2,
\;x_3\in(L_1,L_2),\\
&u_3(x_1,x_2,x_3)\big|_{x_3=L_1}=u_3(x_1,x_2,x_3)\big|_{x_3=L_2}=0,
\;\;\forall(x_1,x_2)\in\R^2.
\end{split}\eeq
Here $\ome\perp u$ and 
\beq\begin{split}\label{cu-ue-ppl}
\ome&=\mbox{curl }u=\nab\times u\\
&=\big(-\patl_1\patl_3\phi(x_1+x_2,x_3),-\patl_2\patl_3
\phi(x_1+x_2,x_3),\{\patl_1^2+\patl_2^2\}\phi(x_1+x_2,x_3)\big).
\end{split}\eeq

Indeed the solution of static incompressible Euler equations
\eqref{SEE1} \eqref{SEE2} with boundary condition \eqref{ue-ppl-b} is
non-unique.
\end{theorem}

This Theorem \ref{Nonu-SEE-ppDb} is proved in Section 2.

\ \ \

Assume that $A(\nab)=e_3\times\nab=(-\patl_2,\patl_1,0)$ and $\psi(t,x)
=0$ in \eqref{Sym-Rep-a2}. Inserting this special \eqref{Sym-Rep-a2}
and \eqref{Def-SRa2-ome} into Navier-Stokes equations \eqref{NS-G-ome},
we obtain

\begin{theorem}[Global Solution of Periodic Pipeline NS Flow]
\label{NSgs-ppDb}  Provided regular function $\phi(t,\xi,\eta)$
satisfies the following equations
\beq\begin{split}\label{phi-ppl}
&\phi_t-\nu(2\patl_{\xi}^2+\patl_{\eta}^2)\phi=0,\\
&\phi(t,\xi+2\pi,\eta)=\phi(t,\xi,\eta),\;\;\forall t\ge0,\;\xi\in\R,
\;\eta\in(L_1,L_2),\\
&\phi(t,\xi,\eta)\big|_{\eta=L_1}=\phi(t,\xi,\eta)\big|_{\eta=L_2}=0,
\;\;\forall t\ge0,\;\xi\in\R,\\
&\phi(t,\xi,\eta)\big|_{t=0}=\phi_0(\xi,\eta).\\
\end{split}\eeq
Let
\beq\begin{split}\label{Def-uns-ppl}
u_{ns}^{pp}(t,x)&=\big(-\patl_2\phi(t,x_1+x_2,x_3),
\patl_1\phi(t,x_1+x_2,x_3),0\big),\\
\end{split}\eeq
then $u=u_{ns}^{pp}$ is the solution of incompressible Navier-Stokes
equations \eqref{NS1} \eqref{NS2} with boundary condition
\beq\label{uns-ppl-b}\begin{split}
&u(t,x_1+2\pi,x_2,x_3)=u(t,x_1,x_2,x_3),\;\;\forall t\ge0,
\;(x_1,x_2)\in\R^2,\;x_3\in(L_1,L_2),\\
&u(t,x_1,x_2+2\pi,x_3)=u(t,x_1,x_2,x_3),\;\;\forall t\ge0,
\;(x_1,x_2)\in\R^2,\;x_3\in(L_1,L_2),\\
&u(t,x_1,x_2,x_3)\big|_{x_3=L_1}=u(t,x_1,x_2,x_3)\big|_{x_3=L_2}=0,
\;\;\forall t\ge0,\;(x_1,x_2)\in\R^2,
\end{split}\eeq
and initial data
\beq\begin{split}\label{uns-ppl-id}
u(t,x)\big|_{t=0}&=\big(-\patl_2\phi_0(x_1+x_2,x_3),
\patl_1\phi_0(x_1+x_2,x_3),0\big).\\
\end{split}\eeq
\end{theorem}

This Theorem \ref{NSgs-ppDb} is easy to prove by using
Theorem \ref{Nonu-SEE-ppDb}. As a consequence of
Theorem \ref{NSgs-ppDb}, we get

\begin{coro}[Typification without Prandtl Layer] \label{No-PPPL}
Under the condition of Theorem \ref{NSgs-ppDb}, we derive that
\[
\lim_{\nu\rightarrow0}u_{ns}^{pp}(t,x)=\big(-\patl_2\phi_0(x_1+x_2,x_3),
\patl_1\phi_0(x_1+x_2,x_3),0\big)=u_{e0}^{pp}(x),
\]
and $u=u_{e0}^{pp}(x)$ is the solution of static incompressible Euler
equations \eqref{SEE1} \eqref{SEE2} with boundary condition
\eqref{ue-ppl-b}.

Here Prandtl layer dose not appear.
\end{coro}

Let us give an interesting exemplification of Theorem \ref{NSgs-ppDb}
by taking
\beq\begin{split}\label{phi0-ppl-f}
\phi_0(\xi,\eta)=&\sum_{j\in J}\sum_{k\in K}\{\al_{jk}\sin(j\xi)+
\beta_{jk}\cos(j\xi)\}\sin\big(\pi k(\eta-L_1)/(L_2-L_1)\big),\\
\end{split}\eeq
where $L_1$, $L_2$, $\al_{jk}$ and $\beta_{jk}$ are constants, $J$ and
$K$ are any finite nonempty subsets of $\{1,2,\cdots\}$. For the 
problem \eqref{phi-ppl}, there exists a unique solution
\beq\begin{split}\label{phit-ppl-f}
\phi(t,\xi,\eta)=\sum_{j\in J}\sum_{k\in K}&e^{-t\nu\{2j^2+\pi^2k^2
/(L_2-L_1)^2\}}\\
&\{\al_{jk}\sin(j\xi)+\beta_{jk}\cos(j\xi)\}
\sin\big(\pi k(\eta-L_1)/(L_2-L_1)\big).\\
\end{split}\eeq
Inserting \eqref{phit-ppl-f} into \eqref{Def-uns-ppl}, we have path
limit
\beq\label{PL-NSs-pplf}\begin{split}
&u_{ef}^{pp\mu}(x)=\lim_{(\nu,t)\in\Gamma_{\mu},(\nu,1/t)\rightarrow
(0,0)}u_{ns}^{pp}(t,x)\\
=&\sum_{j\in J}\sum_{k\in K}e^{-\mu\{2j^2+\pi^2k^2/(L_2-L_1)^2\}}\\
&\Big(-\patl_2\{\al_{jk}\sin\big(j(x_1+x_2)\big)+\beta_{jk}\cos\big(j
(x_1+x_2)\big)\}\sin\frac{\pi k(x_3-L_1)}{L_2-L_1},\\
&\patl_1\{\al_{jk}\sin\big(j(x_1+x_2)\big)+\beta_{jk}\cos\big(j
(x_1+x_2)\big)\}\sin\frac{\pi k(x_3-L_1)}{L_2-L_1},0\Big).\\
\end{split}\eeq
By Theorem \ref{Nonu-SEE-ppDb}, random $u=u_{ef}^{pp\mu}$ is the
solution of static incompressible Euler equations \eqref{SEE1}
\eqref{SEE2} with boundary condition \eqref{ue-ppl-b}. On the other
hand double limit $\lim_{(\nu,1/t)\rightarrow(0,0)}u_{ns}^{pp}(t,x)$
does not exist. These phenomena reveal randomness and turbulence of
incompressible fluids again.

\ \ \

Besides the non-uniqueness of solutions of static Euler equations in
Corollary \ref{Nonu-SEE-pb} and Theorem \ref{Nonu-SEE-ppDb}, the
following non-unique theorem is well known (see \cite{han1,han2} for
instance)

\begin{theorem}[Non-unique Solutions of 2d Euler Equations]
\label{Nonu-SEE-2d} Provided that $r^2=x_1^2+x_2^2$, disc $B_R=
\{(x_1,x_2)\big|r<R\}$, function $\phi(x)=\phi(r)$. For any radial
function $\phi(r)\in C^{\infty}_c(B_R)$, let
\beq\begin{split}\label{Def-ue-2dd}
u_e^{2d}(x_1,x_2)&=\big(-\patl_2\phi(r),\patl_1\phi(r)\big),
\end{split}\eeq
then $u=u_e^{2d}(x_1,x_2)$ is the solution of two dimensional static
incompressible Euler equations with boundary condition
\beq\label{2due-Db}\begin{split}
\sig\cdot u(x_1,x_2)\big|_{(x_1,x_2)\in\patl B_R}=0,
\end{split}\eeq
where $\sig$ is the outer normal of $\patl B_R$.

Indeed the solution of two dimensional static incompressible Euler
equations with boundary condition \eqref{2due-Db} is non-unique too.
\end{theorem}

\section{Proof of Theorem \ref{Nonu-SEE-ppDb}}
\setcounter{equation}{0}

Provided that $A(\nab)=e_3\times\nab=(-\patl_2,\patl_1,0)$,
$\phi(t,x)=\phi(x)$ and $\psi(t,x)=0$ in \eqref{Sym-Rep-a2}.

Theorem \ref{Nonu-SEE-ppDb} can be proved by inserting this special
\eqref{Sym-Rep-a2} \eqref{Def-SRa2-ome} into static Euler equation
\eqref{SEE-G-ome}.

\ \ \

\begin{proof} ~ Inserting $a_1(\nab)=-\patl_2$, $a_2(\nab)=\patl_1$,
$a_3(\nab)=0$ and $\psi=0$ into \eqref{Sym-Rep-a2}
\eqref{Def-SRa2-ome}, we have
\beq\label{uo-SRa2-pp1}\begin{split}
u=&(-\patl_2\phi,\patl_1\phi,0),\\
\ome=&\nab\times u=\big(-\patl_1\patl_3\phi,-\patl_2\patl_3\phi,
(\patl_1^2+\patl_2^2)\phi\big).\\
\end{split}\eeq

It is well-known that the equation \eqref{SEE-G-ome}
is equivalent to
\beq\label{STEE-curl-a2}\begin{split}
&(u\cdot\nab)\ome-(\ome\cdot\nab)u
=\sum_{j=1}^3\patl_j(u^j\ome-\ome^ju)=0.\\
\end{split}\eeq
To verify that $u=u_e$ is the solution of Euler equations \eqref{SEE1}
\eqref{SEE2}, it is necessary to calculate $u^j\ome^k-\ome^ju^k$.
By observation
\beq\label{Ob-uj-ok}\begin{split}
&u^j\ome^k-\ome^ju^k=0,\;\;j=k,\\
&u^j\ome^k-\ome^ju^k=-(u^k\ome^j-\ome^ku^j),
\end{split}\eeq
only three terms $u^2\ome^1-\ome^2u^1$, $u^3\ome^1-\ome^3u^1$ and
$u^3\ome^2-\ome^3u^2$ have to calculate.

Using \eqref{uo-SRa2-pp1}, we derive
\beq\label{uo-21-pp1}\begin{split}
u^2\ome^1-\ome^2u^1
=-\patl_1\phi \patl_1\patl_3\phi-\patl_2\phi \patl_2\patl_3\phi,\\
\end{split}\eeq

\beq\label{uo-31-pp1}\begin{split}
u^3\ome^1-\ome^3u^1
=\patl_2\phi(\patl_1^2\phi+\patl_2^2\phi),\\
\end{split}\eeq

\beq\label{uo-32-pp1}\begin{split}
u^3\ome^2-\ome^3u^2
=-\patl_1\phi(\patl_1^2\phi+\patl_2^2\phi).\\
\end{split}\eeq

Therefore we have
\beq\label{uo1-pp1}\begin{split}
&(u\cdot\nab)\ome^1-(\ome\cdot\nab)u^1=\patl_k(u^k\ome^1-\ome^ku^1)\\
=&\patl_2(u^2\ome^1-\ome^2u^1)+\patl_3(u^3\ome^1-\ome^3u^1)\\
=&\patl_2(-\patl_1\phi \patl_1\patl_3\phi-\patl_2\phi \patl_2\patl_3\phi)
+\patl_3\{\patl_2\phi(\patl_1^2\phi+\patl_2^2\phi)\}\\
=&-\patl_1\patl_2\phi \patl_1\patl_3\phi-\patl_1\phi \patl_1\patl_2\patl_3\phi
+\patl_2\patl_3\phi\patl_1\patl_1\phi+\patl_2\phi\patl_1\patl_1\patl_3\phi,\\
\end{split}\eeq

\beq\label{uo2-pp1}\begin{split}
&(u\cdot\nab)\ome^2-(\ome\cdot\nab)u^2=\patl_k(u^k\ome^2-\ome^ku^2)\\
=&-\patl_1(u^2\ome^1-\ome^2u^1)+\patl_3(u^3\ome^2-\ome^3u^2)\\
=&\patl_1(\patl_1\phi \patl_1\patl_3\phi+\patl_2\phi \patl_2\patl_3\phi)
+\patl_3\{-\patl_1\phi(\patl_1^2\phi+\patl_2^2\phi)\}\\
=&\patl_1\patl_2\phi\patl_2\patl_3\phi+\patl_2\phi\patl_1\patl_2\patl_3\phi
-\patl_1\patl_3\phi\patl_2\patl_2\phi-\patl_1\phi\patl_2\patl_2\patl_3\phi,\\
\end{split}\eeq

\beq\label{uo3-pp1}\begin{split}
&(u\cdot\nab)\ome^3-(\ome\cdot\nab)u^3=\patl_k(u^k\ome^3-\ome^ku^3)\\
=&-\patl_1(u^3\ome^1-\ome^3u^1)-\patl_2(u^3\ome^2-\ome^3u^2)\\
=&-\patl_1\{\patl_2\phi(\patl_1^2\phi+\patl_2^2\phi)\}
+\patl_2\{\patl_1\phi(\patl_1^2\phi+\patl_2^2\phi)\}\\
=&\patl_1\phi\patl_2(\patl_1^2\phi+\patl_2^2\phi)
-\patl_2\phi\patl_1(\patl_1^2\phi+\patl_2^2\phi).\\
\end{split}\eeq

For any smooth function $f(\xi,\eta)$, let $\phi(x_1,x_2,x_3)=
f(x_1+x_2,x_3)$. Then
\beq\label{Pa-phi-f}\begin{split}
\patl_1\phi=\patl_{\xi}f,\;\;\patl_2\phi=\patl_{\xi}f,\;\;
\patl_3\phi=\patl_{\eta}f.\\
\end{split}\eeq
Applying \eqref{Pa-phi-f} in \eqref{uo1-pp1} \eqref{uo2-pp1}
\eqref{uo3-pp1}, we obtain that
\beq\label{uo123-pp1}\begin{split}
&(u\cdot\nab)\ome^1-(\ome\cdot\nab)u^1
=(u\cdot\nab)\ome^2-(\ome\cdot\nab)u^2
=(u\cdot\nab)\ome^3-(\ome\cdot\nab)u^3=0.\\
\end{split}\eeq

Equation \eqref{uo123-pp1} means that $u=u_e^{pp}(x)$ is the solution
of static incompressible Euler equations \eqref{SEE1} \eqref{SEE2}
with boundary condition \eqref{ue-ppl-b}.

\end{proof}


\section*{Acknowledgments}
This work is supported by National Natural Science Foundation of China--NSF,
Grant No.11971068 and No.11971077.


\end{document}